\documentclass[english]{article}
\usepackage[T1]{fontenc}
\usepackage[latin1]{inputenc}
\usepackage{babel}

\usepackage{amsfonts,amsmath,amssymb,theorem}

\font\teneusm=eusm10
\font\seveneusm=eusm7
\font\fiveeusm=eusm5
\newfam\eusmfam
\textfont\eusmfam\teneusm
\scriptfont\eusmfam\seveneusm
\scriptscriptfont\eusmfam\fiveeusm

\renewcommand{\Re}[1]{{\rm Re}(#1)}

\newcommand{\C}{\ensuremath{\mathbb C}}

\newcommand{\R}{\ensuremath{\mathbb R}}

\newcommand{\Z}{\ensuremath{\mathbb Z}}
\newcommand{\N}{\ensuremath{\mathbb N}}

\def\id{\mathop{\text{id}}}

\def\trace{\mathop{\text{trace}}}
\def\Per{\mathop{\text{Per}}}

\newcommand{\cK}{{\cal K}}

\newtheorem{Thm}{Theorem}[section]

\newtheorem{Prop}[Thm]{Proposition}
\newtheorem{Lemma}[Thm]{Lemma}
{\theorembodyfont{\rmfamily}  }
{\theorembodyfont{\rmfamily} \newtheorem{Ex}[Thm]{Example} }

\def\pf{\noindent {\it Proof.}\hskip 8pt}
\def\qed{\hfill{\setlength{\fboxrule}{0.8pt}\setlength{\fboxsep}{1mm}
\fbox{\null}} \vskip 10pt}

\def\3{\ss{}}

\newcommand{\cL}{{\cal L}}
\newcommand{\cF}{{\cal F}}

\newcommand{\PAbl}[2]{\frac{\partial{#1}}{\partial{#2}}}

\def\ssarr{\hbox to 30pt{\rightarrowfill}}
\def\sarr{\hbox to 40pt{\rightarrowfill}}
\def\arr{\hbox to 55pt{\rightarrowfill}}
\def\larr{\hbox to 55pt{\leftarrowfill}}
\def\Arr{\hbox to 80pt{\rightarrowfill}}
\def\mapdown#1{\Big\downarrow\rlap{$\vcenter{\hbox{$\scriptstyle#1$}}$}}
\def\lmapdown#1{\llap{$\vcenter{\hbox{$\scriptstyle#1$}}$}\Big\downarrow}
\def\mapright#1{\smash{\mathop{\arr}\limits^{#1}}}

{}

\begin{document}


\title{The dynamical zeta function and transfer operators for the Kac-Baker model}

\author{J. Hilgert\thanks{
Institut f\"{u}r Mathematik, Technische Universit\"{a}t Clausthal,
38678 Clausthal-Zellerfeld,
Germany. E-mail: hilgert@math.tu-clausthal.de
} \ \ and\  D. Mayer\thanks{
MPI f\"{u}r Mathematik, D-53111 Bonn, on sabbatical from Institut f\"{u}r Theoretische
Physik, Technische Universit\"{a}t Clausthal, 38678 Clausthal-Zellerfeld, Germany.
E-mail: mayer @mpim-bonn.mpg.de or dieter.mayer@tu-clausthal.de
}}

\date{October 30, 2001}

\maketitle
\makeatother

\section{Introduction}

Dynamical zeta functions seem to play a rather special role in the general theory
of zeta functions. They share many properties with the zeta functions of number
theory and algebraic geometry and have lead to new approaches to some of the
outstanding open problems in the theory of these functions. A well known such
problem in this theory is the general Riemann hypothesis for these functions
and a possible spectral interpretation of their zeros and poles in terms of
some dynamical operator attached to the zeta function. There are several examples
like the Selberg or the Artin-Mazur zeta functions where such an operator exists
and which furthermore is indeed related to some dynamical system. In Selberg's
case the underlying dynamical system is the geodesic flow on a surface of constant
negative curvature and the operator in question is the transfer operator of
this flow (see \cite{Ma91}). For the Artin-Weil function the dynamical
system is the Frobenius map on an algebraic variety and the operator is again
the transfer operator for this map (see \cite{Ro86}, \cite{Ru92}).
One should also mention the recent approaches to Riemann's zeta function and
more general $L$-functions by A.\ Connes (\cite{Co96}) and C.\ Deninger
(\cite{De99}) where also a dynamical interpretation of these zeros
is looked for.

In special cases there is even a simple physical interpretation of the zeros
of such a zeta function in terms of energy eigenvalues of a hamiltonian system,
like in the Selberg case, so that it is not surprising that also ideas from
the theory of quantum chaos are being applied now to the Riemann zeros, which
should be connected to the spectrum of some hamiltonian whose classical limit
should be chaotic (cf. \cite{Be86}). One could therefore speculate
if not all the known zeta functions finally turn out to be dynamical zeta functions.

In the same spirit one can interpret also the work of M.\ Gutzwiller in
\cite{Gu82} on the semiclassical quantization of a particle moving in the anisotropic Kepler
potential. There he related the energy spectrum of this particle to the zeros
of the Ruelle zeta function of an abstract dynamical system known in ergodic
theory as a subshift of finite type over two symbols. In statistical mechanics
language this is a lattice spin system of Ising type with a 2-body interaction
of the spins decaying exponentially fast with distance on the lattice. The system
is known in the physical literature as the Kac-Baker model (cf. \cite{Ka59} and
\cite{Ba61}). These authors introduced this model for a better
understanding of phase transitions in weak long-range systems like the van der Waals gas.

It turns out that the Ruelle zeta function of this system, which is just a
generating function for the finite lattice partition functions of the model,
can be expressed in terms of Fredholm determinants of a transfer operator
(see \cite{Ka66}, \cite{Ma80}) and hence its zeros and poles have a nice spectral
interpretation. There exist indeed two completely different such operators for
this system: one is the Ruelle operator  introduced for general
lattice spin systems in dimension one with a rather straightforward physical
interpretation (see \cite{Ru68}), the other one was found by M.\ Kac (see \cite{Ka66}),
using the very special nature of the interaction of this model and its relation
to the Ornstein-Uhlenbeck process. Indeed Kac considered his operator himself
a nice ``trick devoid of any physical significance''.

Whereas Ruelle's operator is acting in a Banach space of observables closely
related to the lattice spin system, the Kac operator is an abstract integral
operator in the Hilbert space \( L^{2}(\mathbb {R},dx) \)  not directly
related to the spin system. It has a rather complicated kernel \( \cal {K}_{\beta } \)
which in the form used also by Gutzwiller in \cite{Gu82} is given by
\begin{equation}
\label{Kacop}
{\cal {K_{\beta }}}(\xi ,\eta )
=\left( \frac{\cosh (\sqrt{\beta J}\xi )\cosh (\sqrt{\beta J}\eta )}
             {\pi \sinh \gamma }\right)^{\frac{1}{2}}
 \exp \left( -\frac{1}{4}\left( \tanh \frac{\gamma }{2}\right)
 (\xi ^{2}+\eta ^{2})-\frac{(\xi -\eta )^{2}}{4\sinh \gamma }\right),
\end{equation}
where \( \beta  \) denotes ``inverse temperature'' and \( J>0 \) and \( 0<\gamma <1 \)
are physical parameters characterizing the strength of the interaction respectively
its decay rate as a function of distance. We should mention that the kernel
used by Kac in \cite{Ka66} is slightly different from the one
given above since he worked with so called fixed boundary conditions whereas
we will use periodic boundary conditions as Gutzwiller did in \cite{Gu82}.

The Ruelle operator \( \cL _{\beta }:B(D)\rightarrow B(D) \) on the other
hand is acting in a Banach space of smooth observables holomorphic in some disc
\( D \) in the complex plane and has the following rather simple looking form
\[
\cL _{\beta }g(z)=e^{\beta Jz}g(\lambda +\lambda z)+e^{-\beta Jz}g(-\lambda +\lambda z),
\]
where \( \lambda =\exp(-\gamma)  \). Various aspects of this operator have been
studied in \cite{Ma80}, \cite{ViMa77} and \cite{PTMT94}. For instance, it defines a family of nuclear
operators holomorphic in the entire complex \( \beta  \)-plane.

The Ruelle zeta function \( \zeta _{R}(z,\beta ) \) of the subshift of finite
type \( (\Omega _{+},\tau ) \), where \( \Omega _{+}=\{1,-1\}^{\mathbb {Z}_{+}} \)
denotes the configuration space over the two symbols \( \left\{ 1,-1\right\}  \)
and \( \tau  \) the shift on the lattice \( \mathbb {Z}_{+} \), is then defined
as
\begin{equation}
\label{zetaR}
\zeta_R (z,\beta )=\exp \sum ^{\infty }_{n=1}\frac{z^{n}}{n}Z_{n}(\beta ),
\end{equation}
where the \( Z_{n}(\beta ) \) denote the so called finite lattice partition
functions of the Kac-Baker model. In his diploma thesis \cite{Mo89}
B. Moritz showed that for real \( \beta  \) this zeta function can be expressed
through Fredholm determinants of both these two operators as
\begin{equation}
\label{Freddet}
\zeta_R (z,\beta )=\frac{\det(1-z\lambda {\cal {G}}_{\beta })}{\det(1-z{\cal {G}}_{\beta })}
                =\frac{\det(1-z\lambda \cL _{\beta })}{\det(1-z\cL _{\beta })},
\end{equation}
where the operator \( {\cal {G_{\beta }}} \) is related to the Kac operator
\( {\cal {K_{\beta }}} \) simply by
\( {\cal {G}}_{\beta }=\frac{1}{\sqrt{\lambda \exp \beta }}{\cal {K}}_{\beta }. \)
{}From this he concluded that the two operators have indeed the same spectra.
However, he could not relate the two operators \( {\cal {G_{\beta }}} \) and
\( \cL _{\beta } \) in an explicit way so that the eigenfunctions of the two
operators could be determined from each other directly.

In the present paper we solve this problem and give explicit formulas relating
the two operators and their eigenfunctions to each other. Hence the two operators
are closely connected to each other and it is not by accident that M.\ Kac could
find his operator. Since the Fredholm determinant of the operator
\( {\cal {L}}_{\beta } \)
is an entire function in the complex variable \( \beta  \) surprisingly also
the Fredholm determinant of the operator \( {\cal {G}}_{\beta } \) can be extended
to an entire function in \( \beta  \). This is certainly not true for the operator
\( {\cal {G}}_{\beta } \) itself.

Our result could throw some light also on a more general problem arising in
the transfer operator method in ergodic theory of dynamical systems: there are
several cases known where the Ruelle operator when acting in some Banach space
of smooth observables turns out, at least numerically, to have real spectrum
for real ``temperature'' \( \beta  \). Unfortunately, there are no general
criteria known in the theory of operators in Banach spaces which would guarantee
such a behaviour. There are very special approaches to this problem for special
systems (see \cite{Rug94}), in several cases however reality of the
spectrum of such transfer operators can be proved by relating the Ruelle operator
to some selfadjoint operator in a Hilbert space. A typical example for this
is the transfer operator for the Gauss map and its generalizations, which when
restricted to a certain Hardy space of holomorphic functions can be conjugated
to an integral operator with symmetric kernel given by Bessel functions (see
\cite{Ma90}, \cite{Is01}). One could therefore speculate
if there does exist always such a relation in cases where the Ruelle operator
has real spectrum for real \( \beta  \). A complete understanding of this rather
mysterious behaviour of transfer operators however is certainly still missing.

In the present example of the Kac-Baker model the Ruelle transfer operator
\( {\cal {L}}_{\beta } \)
acting on a Banach space of holomorphic functions in a disc can be restricted
to a Hilbert space of entire functions in the complex plane square integrable
with respect to the weight function \( \mu _{t}(z)=t\exp(-t\mid z\mid ^{2}) \),
the so called Fock space \( {\cal {H}}L^{2}(\C ,\mu _{t}) \) for
an appropriate value of the parameter \( t \) depending obviously on \( \beta  \).
This Hilbert space is isomorphic to the space \( L^{2}(\R ,d\xi ) \) via the
Segal-Bargmann transform
\( B_{t}:L^{2}(\R ,d\xi )\rightarrow {\cal {H}}L^{2}(\C ,\mu _{t}). \)
It is basically this transformation for \( t=1 \) which relates the Kac-Gutzwiller
operator \( {\cal {G_{\beta }}} \) and the Ruelle operator \( \cL _{\beta } \)
for real \( \beta  \) in a unitary way. For this special parameter value \( t=1 \)
the Segal-Bargmann transformation will be denoted simply by \( B \) and it
is given by
\[
Bf(z)=
2^{\frac{1}{4}}\int _{\R }f(\xi )\exp (2\pi \xi z-\pi \xi ^{2}-\frac{\pi }{2}z^{2})\: d\xi.
\]
This transformation allows us to relate also the eigenfunctions of the two operators
in an explicit way.

By making use of both these two operators we can show that the Ruelle zeta function
\( \zeta_R (1,\beta ) \) has infinitely many zeros and poles on the real axis
unless some unexpected cancellations will take place , whereas there are also
infinitely many trivial zeros at least for the special value \( \lambda =\frac{1}{2} \)
on the line \( \Re \beta =\ln 2 \), again if there are not analogous cancellations.
Indeed, there are no signs for such cancellations numerically. We expect a similar
behaviour also for general \( 0<\lambda <1 \). Obviously it would be interesting
to know if there exists other zeros of this zeta function in the complex
\( \beta  \)--plane
besides the ones mentioned above. If this is not the case, and indeed we have
some numerical hints also for this, this dynamical zeta function would satisfy
some kind of Riemann hypothesis well known for many zeta and $L$-functions of
number theory. This would be a further sign for the basic role dynamical zeta
functions play in the general theory of zeta functions.

\section{The Kac-Gutzwiller and the Ruelle transfer operator}

If \( F=\left\{ 1,-1\right\}  \) denotes the set of possible values of a classical
spin variable \( \sigma  \) the configuration space \( \Omega _{+} \) of all
allowed spin configurations \( \underline{\xi } \) on the half lattice
\( \Z _{+}=\left\{ 0,1,2,\ldots\right\}  \)
is defined as
\( \Omega _{+}=F^{\Z _{+}}=
\left\{ \underline{\xi }=
\left( \xi _{i}\right) _{i\in \Z _{+}}\mid \xi _{i}\in F\: \forall i\in \Z \right\}  \).
The shift \( \tau :\Omega _{+}\rightarrow \Omega _{+} \) is defined as
\( (\tau \underline{\xi })_{i}=\xi _{i+1} \)
for all \( i\in \Z _{+} \) if
\( \underline{\xi }=(\xi _{i})_{i\in \Z _{+}}\in \Omega _{+}. \)
Consider next the two body interaction
\( \Phi (\xi _{i},\xi _{j})=-J\: \xi _{i}\xi _{j}\: \lambda ^{|i-j|} \)
between spins on lattice sites \( i \) and \( j \), where \( J>0 \) determines
the strength of the interaction and \( 0<\lambda <1 \) is its decay rate. Obviously
this interaction decays exponentially fast with distance \(|i-j| \).
The energy \( U_{n}(\xi ) \) of a configuration \( \xi \in \Omega _{+} \) when
restricted to a finite sublattice \( \left[ 0,n-1\right]  \) of length \( n \)
is defined as
\[
U_{n}(\xi )=\sum ^{n-1}_{i=0}\sum ^{\infty }_{j=1}\Phi (\xi _{i},\xi _{i+j}).
\]
A configuration \( \underline{\xi }=(\xi _{i})_{i\in \Z _{+}}\in \Omega _{+} \)
is called periodic with period \( n \) iff \( \xi _{i+n}=\xi _{i} \) for all
\( i\in \Z _{+} \). Denote by \( {\textstyle\Per_{n}} \) the set of periodic configurations
of period \( n \), i.e.
\[
{\textstyle \Per_{n}}=\left\{ \underline{\xi }\in \Omega _{+}\mid\tau ^{n}\underline{\xi }
       =\underline{\xi }\right\} .
\]
The partition functions \( Z_{n}(\beta ) \) for the lattice spin system with
interaction \( \Phi  \) for the finite sublattices \( \left[ 0,n-1\right]  \)
with periodic boundary conditions are defined as
\[
Z_{n}(\beta )=\sum _{\underline{\xi }\in \Per_{n}}
              \exp \left(-\beta U_{n}(\underline{\xi })\right).
\]
Inserting the explicit expression for the interaction \( \Phi  \) for the Kac-Baker
model we get
\[
Z_{n}(\beta )=\sum _{\underline{\xi }\in \Per_{n}}
              \exp \left(\beta J\sum ^{n-1}_{k=0}
                         \sum ^{\infty }_{j=1}\xi _{k}\xi _{k+j}\lambda ^{j}\right).
\]
The physical properties of the spin system are completely determined by these
partition functions. For instance the free energy \( f=f(\beta ) \) is given
by \( f(\beta )=-\beta \: \lim _{n\rightarrow \infty }\frac{1}{n}\ln Z_{n}(\beta ). \)
The main problem of statistical mechanics is then to determine the analytic
properties of this function in the temperature variable \( \beta =\frac{1}{kT} \),
where \( T \) denotes the absolute temperature and \( k \) Boltzmann's constant.
A standard procedure for doing this is the so called transfer matrix method.
There one looks for a matrix whose leading eigenvalue is closely related to
this free energy. Indeed M.\ Kac found for his model an integral operator acting
in some Hilbert space of square integrable functions whose traces can be related
to the partition functions \( Z_{n}(\beta ) \) and hence determine also the
free energy \( f(\beta ) \). As was mentioned before, Kac  in  \cite{Ka66}
did not work with periodic boundary conditions so that the kernel of his integral operator
is slightly different from the operator ${\cal {K_{\beta }}}(\xi ,\eta )$ for periodic
boundary conditions given by (\ref{Kacop}), which we
use in this paper and which was used also by M.\ Gutzwiller in \cite{Gu82}.
The kernel ${\cal {K_{\beta }}}(\xi ,\eta )$
defines a trace class operator in the Hilbert space \( L^{2}(\R ,d\xi ) \).
Obviously \( \cal {K_{\beta }}(\xi ,\eta )=\cal {K_{\beta }}(\eta ,\xi ) \)
and hence the operator \( \cal {K_{\beta }} \) is symmetric for real \( \beta  \).
Its spectrum is therefore real for such \( \beta  \)-values. In the following
we will simply write \( \beta  \) for the expression \( J\beta  \) since \( J \)
is assumed to be fixed. For \( \underline{\xi }\in \Per_{n} \) denote by
\( \underline{\xi }^{p} \)
the configuration \( \underline{\xi }^{p}\in \Omega =F^{\Z } \) on the lattice
\( \Z  \) with \( \xi _{i}^{p}=\xi _{i} \) for \( i\in \Z _{+} \) and
\( \xi _{i+n}^{p}=\xi _{i} \)
for all \( i\in \Z  \). Then one has
\[
\frac{\beta }{2}\sum ^{n-1}_{i=0}\sum ^{+\infty }_{j=-\infty }
                    \xi ^{p}_{i}\xi ^{p}_{j}\exp\left( -\gamma|i-j|\right)
=\frac{1}{2}\underline{\xi }_{n}\cdot \mathbb {A}\underline{\xi }_{n},
\]
where \( \mathbb {A} \) is the \( (n\times n) \)--matrix with matrix elements
\[
\mathbb {A}_{i,j}=\beta \sum ^{+\infty }_{k=-\infty }
\exp(|i-j+nk|),\quad 0\leq i,j\leq n-1
\]
and \( \underline{\xi }_{n}=(\xi _{0},\ldots,\xi _{n-1}). \) On the other hand
the following identity due to H.\ Cram\'er (see \cite{Cr46}, p.\ 118), used also
by M. Kac and M. Gutzwiller, holds :
\[
\exp ({\textstyle\frac{1}{2}}\underline{\xi }_{n}\cdot \mathbb {A}\underline{\xi }_{n})
=(2\pi )^{-\frac{n}{2}}(\det \mathbb {A}^{-1})^{\frac{1}{2}}
 \int ^{+\infty }_{-\infty }dz_{0}\:\ldots\int ^{+\infty }_{-\infty }dz_{n-1}
 \exp(-{\textstyle\frac{1}{2}}\underline{z}\cdot \mathbb {A}^{-1}\underline{z})
 \exp (\underline{\xi }_{n}\cdot \underline{z})
\]
with \( \underline{\xi }_{n}=(\xi _{0},\ldots,\xi _{n-1}) \) and
\( \underline{z}=(z_{0},\ldots,z_{n-1}). \)
{}From this one easily derives the identity
\[
Z_{n}(\beta )=2\sinh\left( \frac{n\gamma }{2}\right)
               \exp\left( -\frac{n\beta }{2}\right)
               \trace\left(\: {\cal {K}_{\beta}}^{n}\right),\]
where \( \cal {K_{\beta }} \) is the integral operator with kernel
\( {\cal {K_{\beta }}}(\xi ,\eta ) \)
as given in the introduction, and where we used the identity
\( \frac{1}{2}(\underline{\xi }_{n}\cdot \mathbb {A}\underline{\xi }_{n})
=-\beta U_{n}(\xi )+\frac{n\beta }{2} \).
Defining next the operator
\[ {\cal {G_{\beta }}}:L^{2}(\R ,d\xi )\rightarrow L^{2}(\R ,d\xi ) \]
with kernel
\( {\cal {G_{\beta }}}(\xi ,\eta )
=\frac{1}{\sqrt{\lambda \exp \beta }}{\cal {K_{\beta }}}(\xi ,\eta ) \)
one finds finally for real \( \beta  \)
\[
Z_{n}(\beta )=(1-\lambda )^{n}\trace\: {\cal {G_{\beta }}}^{n}.
\]
The operator \( {\cal {G_{\beta }}} \) hence serves as a transfer operator
for the Kac-Baker model at least for real \( \beta  \). We call it the Kac-Gutzwiller
operator.

Let us next briefly recall the Ruelle operator for a 1-dimensional lattice spin
system. For this we denote by \( C(\Omega _{+}) \) the space of continuous
observables for the spin system. On this space the following family of linear
operators \( \cL _{\beta }:C(\Omega _{+})\rightarrow C(\Omega _{+}) \) can
be defined
\[
\cL _{\beta }f(\underline{\xi })
=\sum _{\underline{\eta }\in \tau ^{-1}(\underline{\xi })}
   \exp \left(-\beta U_{1}(\underline{\eta })\right)\: f(\underline{\eta }),\]
where \( U_{1}(\underline{\eta }) \) denotes the interaction energy of the
configuration \( \underline{\eta } \) on the sublattice consisting of the lattice
site $0$ only. Inserting the explicit form of \( U_{1}(\underline{\eta }) \)
we find
\[
\cL _{\beta }f(\underline{\xi })
=\sum _{\sigma =+1,-1}
   \exp \left(\beta \sigma \sum ^{\infty }_{i=1}\xi _{i-1}\lambda ^{i}\right)\:
   f((\sigma ,\underline{\xi })),
\]
where we have again replaced \( J\beta  \)  by \( \beta  \).
Obviously the operators \( \cL _{\beta } \) leave invariant the subspace of
functions \( f \) depending on the configuration \( \underline{\xi } \) only
through the variable \( z=\sum ^{\infty }_{i=1}\xi _{i-1}\lambda ^{i} \). Indeed,
they leave invariant also the space of such functions depending holomorphically
on this variable for instance in the disc \( D_{r} \) with
\( r>\frac{\lambda }{1-\lambda }. \)
It was shown in \cite{Ma80} that the operators \( \cL _{\beta } \)
are nuclear for all complex values of \( \beta  \) in the sense of Grothendieck
in the Banach space \( B(D_{r}) \) of functions holomorphic in \( D_{r} \)
and continuous on \( \overline{D_{r}} \) with the supremum norm. Furthermore
the analytic version of the Atiyah-Bott fixed point formula
(see \cite{AtBo67}) shows
that
\[
Z_{n}(\beta )=(1-\lambda ^{n})\trace\cL _{\beta }^{n}.\]
Hence also this Ruelle operator can be used as a transfer operator for the Kac-Baker
model. The Ruelle zeta function $\zeta_R (z,\beta )$
for the Kac-Baker model defined in (\ref{zetaR})
can indeed be expressed through Fredholm determinants of both these two transfer
operators as described in (\ref{Freddet}).
But the Fredholm determinant \( \det (1-z{\cal {L}}_{\beta }) \) of the nuclear
operator \( \cL _{\beta } \) is an entire function both in the variables \( z \)
and \( \beta  \) hence the Ruelle zeta function \( \zeta (z,\beta ) \) is
a meromorphic function in the entire $z$-- and $\beta$--plane. {}From this
it follows immediately (see \cite{Mo89}) that the operators \( \cal {G_{\beta }} \)
and \( \cL _{\beta } \) have the same spectrum for real \( \beta  \). To relate
in this case the two operators and also their eigenfunctions in a more explicit
way we have to investigate the Kac operator in more detail. We should mention
that some of the following arguments have been used already by M. Gutzwiller
in \cite{Gu82}.

\section{Hermite Functions and Mehler's Formula}

Consider the operators \( Z=m_{x}+\frac{1}{2\pi }\PAbl {}{x} \) and
\( Z^{*}=m_{x}-\frac{1}{2\pi }\PAbl {}{x} \)
on \( L^{2}(\R ) \), where \( m_{g} \) is the multiplication operator
\( (m_{g}f)(x)=g(x)f(x) \).
Then the Hermite functions \( h_{k}\in L^{2}(\R ) \) with \( k\in \N _{0} \)
are given by (see \cite{Fo89}, p.51)
\begin{eqnarray*}
h_{0}(x) & = & 2^{\frac{1}{4}}e^{-\pi x^{2}}\\
h_{k}(x) & = & \sqrt{\frac{\pi ^{k}}{k!}}(Z^{*})^{k}h_{0}(x)
              =\frac{2^{\frac{1}{4}}}{\sqrt{k!}}
               \left( \frac{-1}{2\sqrt{\pi }}\right) ^{k}
               e^{-\pi x^{2}}\left( \PAbl {}{x}\right) ^{k}e^{-2\pi x^{2}}.
\end{eqnarray*}
 One has
\[
Z^{*}h_{k}=\sqrt{\frac{k+1}{\pi }}h_{k+1}\quad \mbox {and}\quad
Zh_{k}=\sqrt{\frac{k}{\pi }}h_{k-1}.
\]
The Hermite functions form a Hilbert basis for \( L^{2}(\R ) \) (see \cite{Fo89},
p.53) and the Mehler formula (see \cite{Fo89}, p.55) is the identity
\begin{equation}
\label{Mehler}
\sum _{k=0}^{\infty }\lambda ^{k}h_{k}(x)h_{k}(y)
=\left( \frac{2}{1-\lambda ^{2}}\right) ^{\frac{1}{2}}
\exp \left( \frac{-\pi (1+\lambda ^{2})(x^{2}+y^{2})+4\pi \lambda xy}{1-\lambda ^{2}}\right)
\end{equation}
for \( \lambda \in D_{1} \). Here \( u=\frac{2}{1-\lambda ^{2}} \) lies in
the right half plane, and the sqare root is the branch which is positive for
\( u>0 \).

\begin{Prop}\label{Mehlerprop1}
For \(\gamma\in \{z\in \C\mid \Re z>0\}\) and \(\lambda=e^{-\gamma}\) we have
\[\frac{-\pi(1+\lambda^{2})(x^{2}+y^{2})+4\pi \lambda xy}{1-\lambda^{2}}=
\frac{2\pi}{\sinh\gamma}\left(-\frac{1}{2}(\cosh \gamma)(x^{2}+y^{2})+xy\right).
\]
\end{Prop}

\pf  Set \( a:=\frac{-\pi (1+\lambda ^{2})}{1-\lambda ^{2}} \) and
\( b:=\frac{4\pi \lambda }{1-\lambda ^{2}} \).
Then
\( \frac{b}{a}=\frac{-4e^{-\gamma }}{1+e^{-2\gamma }}
              =\frac{-4}{e^{\gamma }+e^{-\gamma }}
              =\frac{2}{\cosh \gamma } \)
so that
\begin{eqnarray*}
\frac{-\pi (1+\lambda ^{2})(x^{2}+y^{2})+4\pi \lambda xy}{1-\lambda ^{2}}
& = & \frac{4\pi \lambda }{1-\lambda ^{2}}\left( -\frac{1}{2}(\cosh \gamma )
       (x^{2}+y^{2})+xy\right) \\
& = & \frac{4\pi e^{-\gamma }}{1-e^{-2\gamma }}\left( -\frac{1}{2}(\cosh \gamma )
       (x^{2}+y^{2})+xy\right) \\
& = & \frac{2\pi }{\sinh \gamma }\left( -\frac{1}{2}(\cosh \gamma )(x^{2}+y^{2})+xy\right)
\end{eqnarray*}
\qed

\begin{Prop}\label{Mehlerprop2}
\[-\frac{1}{4}\left(\tanh \frac{\gamma}{2}\right)(\xi^{2}+\eta^{2})
           -\frac{(\xi-\eta)^{2}}{4\sinh\gamma}
 =\frac{1 }{2\sinh \gamma}
   \left(-\frac{1}{2}(\cosh\gamma) (\xi^{2}+\eta^{2})+\xi\eta\right).               \]
\end{Prop}

\pf
Using \( (\xi -\eta )^{2}=\xi ^{2}+\eta ^{2}-2\xi \eta  \) we calculate
\begin{eqnarray*}
-\frac{1}{4}\left( \tanh \frac{\gamma }{2}\right) (\xi ^{2}+\eta ^{2})-
 \frac{(\xi -\eta )^{2}}{4\sinh \gamma } & = & -\frac{1}{4}
 \left( \tanh \frac{\gamma }{2}\right) (\xi ^{2}+\eta ^{2})-
 \frac{\xi ^{2}+\eta ^{2}}{4\sinh \gamma }+\frac{1}{2}\frac{\xi \eta }{\sinh \gamma }\\
& = & -\frac{1}{4}\left( \tanh \frac{\gamma }{2}+\frac{1}{\sinh \gamma }\right)
      (\xi ^{2}+\eta ^{2})+\frac{1}{2}\frac{\xi \eta }{\sinh \gamma }\\
& = & -\frac{1}{4}\left( \frac{e^{\gamma }-1}{e^{\gamma }+1}+
      \frac{2}{e^{\gamma }-e^{-\gamma }}\right) (\xi ^{2}+\eta ^{2})+
      \frac{\xi \eta }{e^{\gamma }-e^{-\gamma }}.
\end{eqnarray*}
Set
\( a:=-\frac{1}{4}\left( \frac{e^{\gamma }-1}{e^{\gamma }+1}+
      \frac{2}{e^{\gamma }-e^{-\gamma }}\right)  \)
and \( b:=\frac{1}{e^{\gamma }-e^{-\gamma }} \). Then
\begin{eqnarray*}
\frac{b}{a}
& = & -4\left( \frac{e^{\gamma }-1}{e^{\gamma }+1}+
      \frac{2}{e^{\gamma }-e^{-\gamma }}\right) ^{-1}
      \frac{1}{e^{\gamma }-e^{-\gamma }}\\
 & = & -4\left( \frac{e^{2\gamma }-e^{\gamma }-1+e^{-\gamma }}{e^{\gamma }+1}+2\right) ^{-1}\\
 & = & -4\left( \frac{e^{2\gamma }-e^{\gamma }-1+
        e^{-\gamma }+2e^{\gamma }+2}{e^{\gamma }+1}\right) ^{-1}\\
 & = & -4\left( e^{\gamma }+e^{-\gamma }\right) ^{-1}\\
 & = & -\frac{2}{\cosh \gamma }
\end{eqnarray*}
and therefore
\[
-\frac{1}{4}\left( \tanh \frac{\gamma }{2}\right)
  (\xi ^{2}+\eta ^{2})-\frac{(\xi -\eta )^{2}}{4\sinh \gamma }
=\frac{1}{2\sinh \gamma }\left( -\frac{1}{2}(\cosh \gamma )
 (\xi ^{2}+\eta ^{2})+\xi \eta \right) .
\]
 \qed

\begin{Prop}\label{Mehlerprop3}
 Set \(\lambda=e^{-\gamma}\), \(x=\xi \frac{1}{2\sqrt{\pi}}\)
and     \(y=\eta \frac{1}{2\sqrt{\pi}}\).
Then
\[\sum_{k=0}^{\infty }\lambda^{k}h_{k}(x)h_{k}(y)=
\left(\frac{1}{\sinh \gamma}\right)^{\frac{1}{2}}e^{\frac{\gamma}{2}}
  \exp\left(-\frac{1}{4}\left(\tanh \frac{\gamma}{2}\right)(\xi^{2}+\eta^{2})
           -\frac{(\xi-\eta)^{2}}{4\sinh\gamma}  \right).\]
\end{Prop}

\pf
Combining the Mehler formula (\ref{Mehler}) with Proposition \ref{Mehlerprop1}
and Proposition \ref{Mehlerprop2} we calculate
\begin{eqnarray*}
\sum _{k=0}^{\infty }\lambda ^{k}h_{k}(x)h_{k}(y)
& = & \left( \frac{2}{1-\lambda ^{2}}\right) ^{\frac{1}{2}}
      \exp \left( \frac{2\pi }{\sinh \gamma }
      \left( -\frac{1}{2}(\cosh \gamma )(x^{2}+y^{2})+xy\right) \right) \\
& = & \left( \frac{2}{1-e^{-2\gamma }}\right) ^{\frac{1}{2}}
      \exp \left( \frac{1}{2\sinh \gamma }
      \left( -\frac{1}{2}(\cosh \gamma )(\xi ^{2}+\eta ^{2})+\xi \eta \right) \right) \\
& = & \left( \frac{e^{\gamma }}{e^{\gamma }-e^{-\gamma }}\right) ^{\frac{1}{2}}
      \exp \left( -\frac{1}{4}\left( \tanh \frac{\gamma }{2}\right)
      (\xi ^{2}+\eta ^{2})-\frac{(\xi -\eta )^{2}}{4\sinh \gamma }\right) \\
& = & \left( \frac{1}{\sinh \gamma }\right) ^{\frac{1}{2}}e^{\frac{\gamma }{2}}
      \exp \left( -\frac{1}{4}\left( \tanh \frac{\gamma }{2}\right)
      (\xi ^{2}+\eta ^{2})-\frac{(\xi -\eta )^{2}}{4\sinh \gamma }\right)
\end{eqnarray*}
 \qed

We set
\[
\widetilde{\cK}(\xi ,\eta ):=\left( \frac{1}{\sinh \gamma }\right) ^{\frac{1}{2}}
e^{\frac{\gamma }{2}}\exp \left( -\frac{1}{4}\left( \tanh \frac{\gamma }{2}\right)
(\xi ^{2}+\eta ^{2})-\frac{(\xi -\eta )^{2}}{4\sinh \gamma }\right)
\]
and note that
\begin{equation}
\label{Kacreduced}
\cK_{\beta }(\xi ,\eta )=\left( \frac{\cosh (\sqrt{\beta }\xi )\cosh (\sqrt{\beta }\eta )}
                          {\pi e^{\gamma }}\right) ^{\frac{1}{2}}\widetilde{\cK}(\xi ,\eta ).
\end{equation}

\begin{Lemma}\label{rescale1}
Let \(c\not=0\) and \(a\colon \R\to ]0,\infty[\) be a smooth function. Then
\begin{enumerate}
\item[{\rm(i)}] \(R_{c}\colon L^{2}(\R,a(\frac{\xi}{c})d\xi)\to L^{2}(\R,a(x)dx)\) defined by
               \[(R_{c}f)(x):=\sqrt{|c|} f(cx)\]
               is an isomorphism of Hilbert spaces.
\item[{\rm(ii)}] If an operator \(K\) on \(L^{2}(\R,d\xi)\) is given by an integral
                 kernel \(K(\xi,\eta)\) via
                 \[Kf(\xi)=\int_{\R }K(\xi,\eta) f(\eta)d\eta,\]
                 then the induced operator
                 \(K_{c}:=R_{c}\circ K\circ R_{c}^{-1}\) on \(L^{2}(\R, dx)\)
                 is given by the kernel
                 \[K_{c}(x,y)=cK(cx,cy).\]
\end{enumerate}

\end{Lemma}

\pf

\begin{enumerate}
\item[(i)] This follows from the calculation
\begin{eqnarray*}
\Vert R_{c}f\Vert _{L^{2}(\R ,a(x)dx)}^{2} & = & \int _{\R }|(R_{c})f(x)|^{2}a(x)dx\\
 & = & \int _{\R }{|c|}\, \, |f(cx)|^{2}a(x)dx\\
 & = & \int _{\R }|f(\xi )|^{2}a(\frac{\xi }{c})d\xi \\
 & = & \Vert f\Vert _{L^{2}(\R ,a(\frac{\xi }{c})d\xi )}^{2}.
\end{eqnarray*}

\item[(ii)] This follows from the calculation
\begin{eqnarray*}
\left( (R_{c}\circ K\circ R_{c}^{-1}(f)\right) (x)
& = & \sqrt{|c|}\left( K(R_{c}^{-1}f)\right) (cx)\\
& = & \sqrt{|c|}\int _{\R }K(cx,\eta )(R_{c}^{-1}f)(\eta )d\eta \\
& = & \sqrt{|c|}\int _{\R }K(cx,\eta )\frac{1}{\sqrt{|c|}}f(\frac{1}{c}\eta )d\eta \\
& = & \int _{\R }K(cx,cy)f(y)cdy.
\end{eqnarray*}

\end{enumerate}
\qed

Lemma \ref{rescale1} yields the following commutative diagram:

\[
\begin{array}{ccc}
L^{2}(\R ,d\xi ) & \mapright {K} & L^{2}(\R ,d\xi )\\
\lmapdown {R_{c}} &  & \mapdown {R_{c}}\\
L^{2}(\R ,dx) & \mapright {K_{c}} & L^{2}(\R ,dx)
\end{array}\]

\begin{Ex}\label{Kackernel1} We consider the Kac kernel
\[{\cal{K}}_\beta(\xi,\eta)
=\left(\frac{\cosh(\sqrt{\beta}\xi)\cosh(\sqrt{\beta}\eta)}
      {\pi\sinh \gamma}\right)^{\frac{1}{2}}
 \exp\left(-\frac{1}{4}\left(\tanh \frac{\gamma}{2}\right)(\xi^{2}+\eta^{2})
           -\frac{(\xi-\eta)^{2}}{4\sinh\gamma}\right)
\]
and use the constant \(c=2\sqrt {\pi}\).   Then,
in view of Lemma \ref{rescale1} and  formula (\ref{Kacreduced}),
Proposition \ref{Mehlerprop2} yields
\begin{eqnarray*}
{\cal{K}}_{\beta,c}(x,y)
&=&c{\cal{K}}_\beta(cx,cy)\\
&=&c{\cal{K}}_\beta(\xi,\eta)\\
&=&2\sqrt {\pi}\left(\frac{\cosh(\sqrt{\beta}\xi)\cosh(\sqrt{\beta}\eta)}{\pi e^{\gamma}}
                     \right)^{\frac{1}{2}} \widetilde {K}_\beta(\xi,\eta)\\
&=&2\left(\cosh(2\sqrt {\pi \beta}x)\cosh(2\sqrt {\pi \beta}y)\right)^{\frac{1}{2}}
   \sum_{k=0}^{\infty }e^{-(k+\frac{1}{2})\gamma}h_{k}(x)h_{k}(y)
\end{eqnarray*}
\qed
\end{Ex}

\begin{Lemma}\label{rescale2}
Let \(a\colon \R\to[1,\infty[\) be a smooth function
and \(K\) be a bounded operator on \(L^{2}(\R,dx)\) given by an integral
kernel \(K(x,y)\) via \(Kf(x)=\int_{\R }K(x,y) f(y)dy\). Then
\begin{enumerate}
\item[{\rm(i)}]
The operator \(K \circ m_{\sqrt {a}}\) on \(L^{2}(\R,dx)\)
is an unbounded  integral operator with kernel
\[K(x,y)\sqrt{a(y)}.\]
\item[{\rm(ii)}]
The operator \(m_{\frac{1}{\sqrt {a}}}\circ K\) on \(L^{2}(\R,dx)\)
is a bounded integral operator with kernel
\[\frac{1}{\sqrt{a(x)}}K(x,y).\]
\end{enumerate}
\end{Lemma}

\pf

\begin{enumerate}
\item[(i)] This is immediate from
\[
(K\circ m_{\sqrt{a}}f)(x)=\int _{\R }K(x,y)\sqrt{a(y)}f(y)dy\]
 and \( \Vert \sqrt{a}f\Vert _{L^{2}(\R ,dx)}=\Vert f\Vert _{L^{2}(\R ,a(x)dx)}<\infty  \).
\item[(ii)] Here we calculate
\begin{eqnarray*}
(m_{\frac{1}{\sqrt{a}}}\circ Kf)(x) & = & \frac{1}{\sqrt{a(x)}}(Kf)(x)\\
 & = & \frac{1}{\sqrt{a(x)}}\int _{\R }K(x,y)f(y)dy\\
 & = & \int _{\R }\frac{1}{\sqrt{a(x)}}K(x,y)f(y)dy
\end{eqnarray*}

\end{enumerate}
\qed

For a smooth map \( a\colon \R \rightarrow [1,\infty [ \) Lemma \ref{rescale2}
yields the following commutative diagram
\[
\begin{array}{ccc}
L^{2}(\R ,a(x)dx) & \mapright {K} & L^{2}(\R ,dx)\\
\lmapdown {m_{\sqrt{a}}} &  & \mapdown {\id }\\
L^{2}(\R ,dx) & \mapright {K'} & L^{2}(\R ,dx)\\
\lmapdown {\id } &  & \mapdown {m_{\frac{1}{\sqrt{a}}}}\\
L^{2}(\R ,dx) & \mapright {K''} & L^{2}(\R ,dx)
\end{array}\]
 with integral operators \( K' \) and \( K'' \) given by the kernels
 \( K'(x,y)=K(x,y)\frac{1}{\sqrt{a(y)}} \)
and
\begin{equation}
\label{rescaleEq}
K''(x,y)=\frac{1}{\sqrt{a(x)}}K(x,y)\frac{1}{\sqrt{a(y)}}.
\end{equation}

\begin{Ex}\label{Kackernel2} We consider the rescaled Kac kernel
\[\cK_{\beta,c}(x,y)
=2\left(\cosh(2\sqrt {\pi\beta }x)\cosh(2\sqrt {\pi\beta }y)\right)^{\frac{1}{2}}
   \sum_{k=0}^{\infty }e^{-(k+\frac{1}{2})\gamma}h_{k}(x)h_{k}(y)\]
with \(c=2\sqrt {\pi}\) from Example \ref{Kackernel1}
and use the function \(a(x)=\cosh (2\sqrt{\pi\beta }x)\).
Then (\ref{rescaleEq}) shows
\[\cK_{c}''(x,y)
  =2\sum_{k=0}^{\infty }e^{-(k+\frac{1}{2})\gamma}h_{k}(x)h_{k}(y),\]
where we have omitted the index \(\beta\) since the kernel does not depend on it anymore. Thus the Hilbert basis \((h_{k})_{k\in \N_{0}}\)
diagonalizes
\(\cK_{c}''\colon L^{2}(\R,dx)\to L^{2}(\R,dx)\) and the corresponding eigenvalues
are given by
\[\cK_{c}'' h_{k}=  2 e^{-(k+\frac{1}{2})\gamma} h_{k}.
\]
Lemma \ref{rescale1} now yields the following commutative diagram
\[
\begin{array}{ccc}
L^{2}(\R,a(\frac{\xi}{c})d\xi)&\mapright{\cK_\beta}&L^{2}(\R,d\xi)\\
\lmapdown{R_{c}}&&\mapdown{R_{c}}\\
L^{2}(\R,a(x)dx)&\mapright{\cK_{\beta,c}}&L^{2}(\R,dx)\\
\lmapdown{m_{\sqrt {a}}}&&\mapdown{\id}\\
L^{2}(\R,dx)&\mapright{\cK_{\beta,c}'}&L^{2}(\R,dx)\\
\lmapdown{\id}&&\mapdown{m_{\frac{1}{\sqrt {a}}}}\\
L^{2}(\R,dx)&\mapright{\cK_{c}''}&L^{2}(\R,dx)
\end{array}
\]

Note that \(e^{\pi x^{2}} h_{k}(x)\) is a polynomial of degree \(k\) in \(x\)
(see \cite{Fo89}, p.52). In view of the estimate
\[\sqrt{\cosh R\xi}  \le e^{\frac{1}{2}\left|R\xi \right| }\]
this shows that the functions
\(\xi\mapsto \sqrt{\cosh\sqrt{\beta} \xi}\ h_{k}(\xi)\) are in \(L^{2}(\R,d\xi)\).
Therefore all the \(h_{k}\) are contained in
\(L^{2}(\R,\cosh (\sqrt{\beta}\xi) d\xi)\).
\qed
\end{Ex}

\begin{Prop}\label{commuteProp} For \(c\in\R\setminus\{0\}\)
und \(s\in \R\) we have
\[m_{\cosh( sx)}\circ R_{c} =R_{c}\circ m_{\cosh (\frac{ s}{c}\xi)}.\]
\end{Prop}

\pf

\begin{eqnarray*}
(m_{\cosh (sx)}\circ R_{c}f)(x) & = & \cosh (sx)\sqrt{|c|}f(cx)\\
 & = & \sqrt{|c|}\cosh (\frac{s}{c}cx)f(cx)\\
 & = & (R_{c}\circ m_{\cosh (\frac{s}{c}\xi )}f)(x)
\end{eqnarray*}
 \qed

\begin{Prop}\label{rescale3} Set \(g(\xi)=\cosh(R'\xi)^{t}\) for \(t\in \R^{+}\) and
\(R'\in \R\). Then the Kac operator
\(\cK_\beta\colon L^{2}(\R,d\xi)\to L^{2}(\R,d\xi)\) satisfies
\(\cK_\beta\big(L^{2}(\R,d\xi)\big)\subseteq L^{2}(\R,g(\xi)d\xi)\).
\end{Prop}

\pf
The estimate
\[
|g(x)|\leq e^{t|R'\xi |}\]
 shows that there exist positive constants \( c_{1} \) and \( c_{2} \) such
that
\[
|g(\xi )\cK_\beta(\xi ,\eta )|\leq c_{1}e^{-c_{2}(\xi ^{2}+\eta ^{2})}.\]
 Therefore there exists a positive constant \( C \) such that
\[
\int _{\R }|\cK_\beta(\xi ,\eta )|g(\xi )d\xi <C\quad \forall \eta \in \R \]
 \(  \)and
\[
\int _{\R }|\cK_\beta(\xi ,\eta )|d\eta <C\quad \forall \xi \in \R .\]
 Now \cite{Fo84}, Theorem 6.18, implies that \( \cK_\beta(\xi ,\eta ) \) defines a
bounded integral operator
\[
\cK_\beta\colon L^{2}(\R ,d\eta )\rightarrow L^{2}(\R ,g(\xi )d\xi )\]
 with operator norm bounded by \( C \). In particular \( \cK_\beta f\in L^{2}(\R ,g(\xi )d\xi ) \).
\qed

\section{Fock Space and Segal-Bargmann Transformation}

\label{Bargmann}

For \( t>0 \) the Fock space \( {\cal {H}}L^{2}(\C ,\mu _{t}) \)
consists of all entire functions \( F\colon \C \rightarrow \C  \) such that
\[
\Vert F\Vert _{t}^{2}:=\int _{\C }|F(z)|^{2}\mu _{t}(z)dz<\infty \]
 where \( \mu _{t} \) denotes the weight \( \mu _{t}=t\exp -\pi t\mid z\mid ^{2} \).
The space \( {\cal {H}}L^{2}(\C ,\mu _{t}) \) with this norm is a Hilbert space
(see \( \left[ Fo89\right]  \), p.46 ). Moreover the polynomials
\[
\zeta ^{t}_{k}(z):=\sqrt{\frac{\pi ^{k}}{t^{k}k!}}z^{k}\quad \forall k\in \N _{0}\]
form an orthonormal basis for the space \( {\cal {H}}L^{2}(\C ,\mu _{t}) \).
Of special interest for the following discussion is the case \( t=1. \) In
this case the space \( {\cal {H}}L^{2}(\C ,\mu _{1}) \) will be denoted simply
by \( \cF  \) and we call it the Fock space (see \cite{Fo89}, p.40
):
\[
\cF =\left\{ F:\C \rightarrow \C \mid  \: \mbox{F entire with}\:
                \int _{\C }|F(z)|^{2}\exp -(\pi z^{2})\: dz<\infty \right\}
\]
Its basis \( \zeta _{k}^{1} \) will be denoted simply by \( \zeta _{k} \) with
$k=0,1,\ldots$ .

There is an isomorphism
\( B_{t}\colon L^{2}(\R ,d\xi )\rightarrow {\cal {H}}L^{2}(\C ,\mu _{t}) \),
the so called Segal-Bargmann transform, defined by
\[
B_{t}f(z):=\left({\textstyle\frac{2}{t}}\right)^{\frac{1}{4}}\int _{\R }f(\xi )
           \exp (2\pi \xi z-{\textstyle \frac{\pi }{t}\xi ^{2}-\frac{\pi t}{2}z^{2}})d\xi .
\]
It is a Hilbert space isomorphism such that \( B_{t}h^{t}_{k}=\zeta ^{t}_{k} \)
(see \cite{Fo89}, p.40 and p.51), where the Hermite functions \( h^{t}_{k} \)
are generalizations of the Hermite functions \( h^{1}_{k}=h_{k} \).

Recall the operator \( \cK_{c}'' \) on \( L^{2}(\R ,dx) \) from Example \ref{Kackernel2}.
Obviously the Segal-Bargmann transform \( B:=B_{1} \) induces a linear operator in
the space \( \cF  \). Indeed one finds

\begin{Prop} For \(\lambda=e^{-\gamma}\) the operator
\(M_{\lambda}:=B\circ \cK_{c}''\circ B^{-1}\colon  \cF\to\cF\)
is given by
\[(M_{\lambda }F)(z)=2\sqrt{\lambda }\ F(\lambda z).\]
\end{Prop}

\pf
Using Example \ref{Kackernel2} and the homogeneity of \( \zeta _{k} \) we
calculate
\begin{eqnarray*}
(M_{\lambda }\zeta _{k})(z) & = & B(\cK_{c}''h_{k})(z)\\
 & = & B(2e^{-(k+\frac{1}{2})\gamma }h_{k})(z)\\
 & = & 2e^{-(k+\frac{1}{2})\gamma }\zeta _{k}(z)\\
 & = & 2e^{-\frac{1}{2}\gamma }\zeta _{k}(e^{-\gamma }z).
\end{eqnarray*}
 Since \( (\zeta _{k})_{k\in \N _{0}} \) is a Hilbert basis for \( \cF  \)
this implies the claim. \qed

The functions \( \zeta _{k}(z), k=0,1,\ldots \) are the eigenfunctions of the
operator \( M_{\lambda } \) with eigenvalue \( \rho _{k}=2\lambda ^{k+\frac{1}{2}} \).
Since this operator is nuclear its spectrum consists of these numbers. The
above proposition shows that the following diagram commutes:

\[
\begin{array}{ccc}
L^{2}(\R ,dx) & \mapright {\cK_{c}''} & L^{2}(\R ,dx)\\
\mapdown {B} &  & \mapdown {B}\\
\cF  & \mapright {M_{\lambda }} & \cF
\end{array}\]

\begin{Prop} Let \(\tau_{r}\colon L^{2}(\R,dx)\to L^{2}(\R,dx)\) for \(r\in \R\)
be the translation defined by
\[\tau_{r}f(x)=f(x-r),\]
and \(\mu_{r}\colon L^{2}(\R,dx)\to L^{2}(\R,dx)\) for \(r\in\R\) be the multiplication
defined by
\[\mu_{r}f(x)=r f(x).\]
Then we have
\[(Z^{*})^{k}\circ \tau_{r}=\tau_{r}\circ(Z^{*}+\mu_{r})^{k}\quad \forall k\in \N_{0}.\]
In particular we find
\[(Z^{*})^{k}\circ \tau_{r} h_{0}
=\sum_{l=0}^{k} \begin{pmatrix}k\\ l\end{pmatrix} r^{l}\tau_{r} \tilde h_{k-l}.\]
\end{Prop}

\pf
For the first claim it suffices to prove the case \( k=1 \). So we calculate

\begin{eqnarray*}
\left( (Z^{*}\circ \tau _{r})f\right) (x)
& = & x\tau _{r}f(x)-\frac{1}{2\pi }(\tau _{r})'(x)\\
& = & xf(x-r)-\frac{1}{2\pi }f'(x-r)\\
& = & (x-r)f(x-r)-\frac{1}{2\pi }f'(x-r)+rf(x-r)\\
& = & (\tau _{r}\circ Z^{*}f)(x)+(\mu _{r}\circ \tau _{r}f)(x)\\
& = & (\tau _{r}\circ Z^{*}f)(x)+(\tau _{r}\circ \mu _{r}f)(x)\\
& = & \left( \tau _{r}\circ (Z^{*}+\mu _{r})f\right) (x).
\end{eqnarray*}
 Now we calculate
\begin{eqnarray*}
(Z^{*})^{k}\circ \tau _{r}h_{0} & = & \tau _{r}\circ (Z^{*}+\mu _{r})^{k}h_{0}\\
& = & \tau _{r}\circ \sum _{l=0}^{k}\left( _{l}^{k}\right) \mu ^{l}_{r}(Z^{*})^{k-l}h_{0}\\
& = &\tau _{r}\circ \sum _{l=0}^{k}\left( _{l}^{k}\right) r^{l}\widetilde{h}_{k-l}\\
& = &\sum _{l=0}^{k}\left( _{l}^{k}\right) r^{l}\: \tau _{r}\widetilde{h}_{k-l}
\end{eqnarray*}
\qed

\begin{Prop}\label{multexp1}
For each \(k\in \N_{0}\) we have
\[(Z^{*})^{k}\circ m_{e^{sx}}
  =m_{e^{sx}}\circ (Z^{*}-\textstyle{\frac{s}{2\pi}})^{k}\]
and
\[m_{e^{sx}}\circ  (Z^{*})^{k}
  = (Z^{*}+\textstyle{\frac{s}{2\pi}})^{k} \circ m_{e^{sx}}.\]
\end{Prop}

\pf
It suffices to prove the first equality in the case \( k=1 \). So we calculate
\begin{eqnarray*}
Z^{*}(e^{sx}f)(x) & = & xe^{sx}f(x)-\frac{1}{2\pi }(e^{sx}f(x))'\\
& = & xe^{sx}f(x)-\frac{1}{2\pi }(se^{sx}f(x)+e^{sx}f'(x))\\
& = & e^{sx}(Z^{*}f)(x)-\frac{s}{2\pi }e^{sx}f(x)\\
& = & m_{e^{sx}}\left( Z^{*}-\frac{s}{2\pi }\right) f(x)
\end{eqnarray*}
 \qed

Set
\[
\widetilde{h}_{k}:=(Z^{*})^{k}h_{0}.\]

\begin{Prop}\label{multexp2} For \(s\in \R\) we have
\[m_{e^{sx}} \widetilde {h}_{k}
=e^{\frac{s^{2}}{4\pi}}
    \sum_{j=0}^{k}  \begin{pmatrix}k\\ j\end{pmatrix}
                             \left(\frac{s}{\pi}\right)^{k-j}
                 \tau_{\frac{s}{2\pi}} \tilde {h}_{j}.
\]
\end{Prop}

\pf

\begin{eqnarray*}
m_{e^{sx}}h_{0}(x) & = & e^{sx}2^{\frac{1}{4}}e^{-\pi x^{2}}\\
 & = & 2^{\frac{1}{4}}e^{sx-\pi x^{2}}\\
 & = & 2^{\frac{1}{4}}e^{-\pi (\frac{s}{2\pi }-x)^{2}+\left( \frac{s}{2\pi }\right) ^{2}\pi }\\
 & = & e^{\left( \frac{s}{2\pi }\right) ^{2}\pi }2^{\frac{1}{4}}e^{-\pi (\frac{s}{2\pi }-x)^{2}}\\
 & = & e^{\left( \frac{s}{2\pi }\right) ^{2}\pi }h_{0}(x-\frac{s}{2\pi })\\
 & = & e^{\frac{s^{2}}{4\pi }}\tau _{\frac{s}{2\pi }}h_{0}(x).
\end{eqnarray*}
 Now we use Proposition \ref{multexp1} to calculate

\begin{eqnarray*}
m_{e^{sx}}\widetilde{h}_{k} & = & m_{e^{sx}}(Z^{*})^{k}h_{0}\\
 & = & (Z^{*}+\frac{s}{2\pi })^{k}m_{e^{sx}}h_{0}\\
 & = & e^{\frac{s^{2}}{4\pi }}(Z^{*}+\frac{s}{2\pi })^{k}\tau _{\frac{s}{2\pi }}h_{0}\\
 & = & e^{\frac{s^{2}}{4\pi }}\sum ^{k}_{l=0}\left( _{l}^{k}\right)
       (\frac{s}{2\pi })^{k-l}\: (Z^{*})^{l}\: \tau _{\frac{s}{2\pi }}h_{0}\\
 & = & e^{\frac{s^{2}}{4\pi }}\sum ^{k}_{l=0}\left( _{l}^{k}\right)
       (\frac{s}{2\pi })^{k-l}\: \tau _{\frac{s}{2\pi }}
       (Z^{*}+\mu _{\frac{s}{2\pi }})^{l}h_{0}\\
 & = & e^{\frac{s^{2}}{4\pi }}\sum ^{k}_{l=0}\left( _{l}^{k}\right)
       (\frac{s}{2\pi })^{k-l}\: \sum ^{l}_{j=0}\left( _{j}^{l}\right)
       (\frac{s}{2\pi })^{l-j}\tau _{\frac{s}{2\pi }}(Z^{*})^{j}h_{0}\\
 & = & e^{\frac{s^{2}}{4\pi }}\sum ^{k}_{l=0}\sum ^{l}_{j=0}\left( _{l}^{k}\right)
       \left( _{j}^{l}\right) (\frac{s}{2\pi })^{k-j}\: \tau _{\frac{s}{2\pi }}
       \widetilde{h}_{j}\\
 & = & e^{\frac{s^{2}}{4\pi }}\sum ^{k}_{j=0}\sum ^{k}_{l=j}\left( _{l}^{k}\right)
       \left( _{j}^{l}\right) (\frac{s}{2\pi })^{k-j}\: \tau _{\frac{s}{2\pi }}
       \widetilde{h}_{j}
\end{eqnarray*}

On the other hand
\[
\sum ^{k}_{l=j}\left( _{l}^{k}\right) \left( _{j}^{l}\right) \: =\sum ^{k-j}_{n=0}
  \left( _{n+j}^{k}\right) \left( _{j}^{n+j}\right)
=\sum ^{k-j}_{n=0}\left( _{n}^{k-j}\right) \left( _{j}^{k}\right)
=2^{k-j}\left( _{j}^{k}\right)
\]
which inserted in the calculation above gives the claim.
\qed

\begin{Prop} \label{Bargtrans} For \(r\in \R\) we have
\[B\circ \tau_{r}=m_{e^{\pi rz-\frac{\pi}{2} r^{2}}}\circ \tau_{r}\circ B.\]
Here \(\tau_{r}\colon \cF\to\cF\) is also defined by \(\tau_{r}F(z)=F(z-r)\).
In other words, for \(F\in \cF\) we have
\[(B\circ \tau_{r}\circ B^{-1})F(z)= e^{-\frac{\pi}{2}r^{2}} e^{\pi rz}F(z-r).\]
\end{Prop}

\pf

\begin{eqnarray*}
\left( (B\circ \tau _{r})f\right) (z)
& = & 2^{\frac{1}{4}}\int _{\R }f(\xi -r)\exp (2\pi \xi z-\pi \xi ^{2}-\frac{\pi }{2}z^{2})d\xi \\
& = & 2^{\frac{1}{4}}\int _{\R }f(\xi )\exp (2\pi (\xi +r)-\pi (\xi +r)^{2}-\frac{\pi }{2}z^{2})d\xi \\
& = & 2^{\frac{1}{4}}\int _{\R }f(\xi )\exp (2\pi \xi (z-r)-\pi \xi ^{2}-\frac{\pi }{2}(z-r)^{2}+\pi rz-\frac{\pi }{2}r^{2})d\xi \\
& = & e^{\pi rz-\frac{\pi }{2}r^{2}}2^{\frac{1}{4}}\int _{\R }f(\xi )\exp (2\pi \xi (z-r)-\pi \xi ^{2}-\frac{\pi }{2}(z-r)^{2})d\xi \\
& = & e^{\pi rz-\frac{\pi }{2}r^{2}}Bf(z-r)\\
& = & e^{\pi rz-\frac{\pi }{2}r^{2}}\tau _{r}Bf(z)
\end{eqnarray*}
 \qed

\begin{Prop} \label{multexp3}For \(s\in \R\) we have
\[B\circ m_{e^{sx}}\circ B^{-1}
= e^{\frac{1}{8\pi}s^{2}} m_{e^{\frac{sz}{2}}}\circ \tau_{-s}.\]
In other words, for \(F\in \cF\) and \(z\in \C\) we have
\[(B\circ m_{e^{sx}}\circ B^{-1} F)(z)=
 e^{\frac{s^{2}}{8\pi}} e^{\frac{sz}{2}} F(z+\frac{s}{2\pi}).\]
\end{Prop}

\pf
Using Proposition \ref{multexp2}, Proposition \ref{Bargtrans} and the fact
that \( (B\widetilde{h}_{j})(z)=z^{j} \) we calculate
\begin{eqnarray*}
B\circ m_{\exp sx} \widetilde{h_{k}}
&=&\exp \left(\frac{s^{2}}{4\pi }\right)\sum ^{k}_{j=0}\left( _{j}^{k}\right)
    \left(\frac{s}{\pi }\right)^{k-j}    B\circ \tau _{\frac{s}{2\pi }}\widetilde{h_{j}}\\
&=&\exp \left(\frac{s^{2}}{4\pi }\right)\sum ^{k}_{j=0}\left( _{j}^{k}\right)
     \left(\frac{s}{\pi }\right)^{k-j} \exp \left(\frac{sz}{2}-\frac{s^{2}}{8\pi }\right)\!
     \tau _{\frac{s}{2\pi }}B\! \widetilde{h_{j}}\\
&=&\exp \left(\frac{s^{2}}{8\pi }\right)\exp \left(\frac{sz}{2}\right)
    \: \sum ^{k}_{j=0}\left( _{j}^{k}\right)
   (\frac{s}{\pi })^{k-j}\: \tau _{\frac{s}{2\pi }}z^{j}\\
&=&\exp \left(\frac{s^{2}}{8\pi }\right)\exp \left(\frac{sz}{2}\right)\:
    \sum ^{k}_{j=0}\left( _{j}^{k}\right) (\frac{s}{\pi })^{k-j}(z-\frac{s}{2\pi })^{j}\\
&=&\exp \left(\frac{s^{2}}{8\pi }\right)\exp \left(\frac{sz}{2}\right)\:
   (z+\frac{s}{2\pi })^{k}
\end{eqnarray*}
and since the \( \widetilde{h_{k}} \) form a basis of \( \cF  \) this implies
the claim.
\qed

Consider next the unbounded operator \( C_{s}\colon \cF \rightarrow \cF  \)
defined as
\[
C_{s}:=B\circ m_{\cosh (sx)}\circ B^{-1}.\]
Then one finds
\begin{equation}
\label{Csdef}
\textstyle C_{s}F(z)={\frac{1}{2}}e^{\frac{s^{2}}{8\pi }}\left( e^{\frac{sz}{2}}F(z+\frac{s}{2\pi })+e^{-\frac{sz}{2}}F(z-\frac{s}{2\pi })\right)
\end{equation}
 and we have the following commutative diagram
\[
\begin{array}{ccc}
L^{2}(\R ,dx) & \mapright {m_{\cosh (sx)}} & L^{2}(\R ,dx)\\
\mapdown {B} &  & \mapdown {B}\\
\cF  & \mapright {C_{s}} & \cF
\end{array}\]

\begin{Prop} \label{CsProp}
For \(s\in \R\) and \(\lambda=e^{-\gamma}\) we have
\[(C_{s}\circ M_{\lambda }F)(z)
=\sqrt{\lambda }e^{\frac{1}{8\pi^{2}}s^{2}}\left(e^{\frac{sz}{2}}F(\lambda z+\lambda \frac{s}{2\pi})
                                      +e^{-\frac{sz}{2}}F(\lambda z- \lambda \frac{s}{2\pi})\right).\]
\end{Prop}

\pf

\begin{eqnarray*}
\textstyle (C_{s}\circ M_{\lambda }F)(z) & = & {\frac{1}{2}}e^{\frac{s^{2}}{8\pi }}\left( e^{\frac{sz}{2}}(M_{\lambda }F)(z+\frac{s}{2\pi })+e^{-\frac{sz}{2}}(M_{\lambda }F)(z-\frac{s}{2\pi })\right) \\
 & = & \sqrt{\lambda }e^{\frac{s^{2}}{8\pi }}\left( e^{\frac{sz}{2}}F(\lambda z+\lambda \frac{s}{2\pi })+e^{-\frac{sz}{2}}F(\lambda z-\lambda \frac{s}{2\pi })\right)
\end{eqnarray*}
 \qed

For \( \alpha \neq 0 \) consider the map
\[
\nu _{\alpha }:{\cal {H}}L^{2}(\C ,\mu _{\frac{1}{ \mid \alpha \mid ^{2}}})\rightarrow \cF ,\]
defined as
\[
(\nu _{\alpha }F)(z):=F(\alpha z).\]
Then one finds for the operator
\[
C_{s}\circ M_{\lambda }\circ \nu _{\alpha }:{\cal {H}}L^{2}(\C ,\mu _{\frac{1}{\mid \alpha \mid ^{2}}})\rightarrow \cF \]

\begin{eqnarray*}
(C_{s}\circ M_{\lambda }\circ \nu _{\alpha }F)(z)
&=&\sqrt{\lambda }e^{\frac{s^{2}}{8\pi }}\left( e^{\frac{sz}{2}}
   (M_{\lambda }\circ\nu _{\alpha }\circ F)(z+\frac{s}{2\pi })
    +e^{-\frac{sz}{2}}(M_{\lambda }\circ\nu _{\alpha }\circ F)(z-\frac{s}{2\pi })\right)\\
&=&\sqrt{\lambda }e^{\frac{s^{2}}{8\pi }}\left( e^{\frac{sz}{2}}
    F(\alpha \lambda z+\alpha \lambda \frac{s}{2\pi })
    +e^{-\frac{sz}{2}}F(\alpha \lambda z-\alpha \lambda \frac{s}{2\pi })\right)
\end{eqnarray*}
respectively for the operator
\[
\nu _{\alpha ^{-1}}\circ C_{s}\circ M_{\lambda }\circ \nu _{\alpha }:{\cal {H}}L^{2}(\C ,\mu _{\frac{1}{\mid \alpha \mid ^{2}}})\rightarrow {\cal {H}}L^{2}(\C ,\mu _{\frac{1}{\mid \alpha \mid ^{2}}})\]

\[
(\nu _{\alpha ^{-1}}\circ C_{s}\circ M_{\lambda }\circ \nu _{\alpha }F)(z)=\sqrt{\lambda }e^{\frac{s^{2}}{8\pi }}\left( e^{\frac{sz}{2\alpha }}F(\lambda z+\alpha \lambda \frac{s}{2\pi })+e^{-\frac{sz}{2\alpha }}F(\lambda z-\alpha \lambda \frac{s}{2\pi })\right) \]
Next chose the parameters \( s \) and \( \alpha  \) as follows:
\[
s=2\sqrt{\pi \beta }\: \: \mbox{ and }\: \: \alpha =\sqrt{\frac{\pi }{\beta }.}\]
Obviously, \( \alpha  s=2\pi  \) and therefore one gets for these parameters
for the operator
\[
\nu _{\sqrt{\frac{\beta }{\pi }}}\circ C_{2\sqrt{\beta \pi }}\circ M_{\lambda }\circ
\nu _{\sqrt{\frac{\pi }{\beta }}}:
{\cal {H}}L^{2}(\C ,\mu _{\frac{\mid \beta \mid }{\pi }})\rightarrow
{\cal {H}}L^{2}(\C ,\mu _{\frac{\mid \beta \mid }{\pi }})\]
defined by
\[
(\nu _{\sqrt{\frac{\beta }{\pi }}}\circ C_{2\sqrt{\beta \pi }}\circ M_{\lambda }\circ
\nu _{\sqrt{\frac{\pi }{\beta }}}F)(z)=\sqrt{\lambda }e^{\frac{\beta }{2}}
\left( e^{\beta z}F(\lambda z+\lambda )+e^{-\beta z}F(\lambda z-\lambda )\right) .\]

\begin{Prop} \label{Ruelle}
For real \(\beta\) and $c=2\sqrt \pi$ the operators
\[
\cL _{\beta }:{\cal {H}}L^{2}(\C ,\mu _{\frac{\mid \beta \mid }{\pi }})\rightarrow
{\cal {H}}L^{2}(\C ,\mu _{\frac{\mid \beta \mid }{\pi }})
\]
and
\[
{\textstyle \frac{1}{\sqrt{\lambda \exp \beta }}}
m_{\cosh (2\sqrt{\beta \pi }x)}\circ \cK^{''}_{c}:L^{2}(\R ,\, dx)\rightarrow
L^{2}(\R ,\, dx)
\]
are conjugate operators.
\end{Prop}

\pf
Indeed one finds after inserting the definitions of the two operators
\( C_{2\sqrt{\beta \pi }} \)
and \( M_{\lambda } \)
\[
C_{2\sqrt{\beta \pi }}\circ M_{\lambda }
=B\circ m_{\cosh 2\sqrt{\beta \pi }x}\circ B^{-1}\circ B\circ \cK_{c}^{''}\circ B^{-1}
=B\circ m_{\cosh 2\sqrt{\beta \pi }x}\circ \cK_{c}^{''}\circ B^{-1}
\]
which shows that \( C_{2\sqrt{\beta \pi }}\circ M_{\lambda } \) is conjugate
to the operator \( m_{\cosh (2\sqrt{\beta \pi }x)}\circ \cK^{''}_{c} \). Hence
\( f\in L^{2}(\R ,\, dx) \) is an eigenfunction of the operator
\( m_{\cosh (2\sqrt{\beta \pi }x)}\circ \cK^{''}_{c} \)
with eigenvalue \( \rho  \) iff the function \( Bf \) is an eigenfunction
of the operator \( C_{2\sqrt{\beta \pi }}\circ M_{\lambda }:\cF \rightarrow \cF  \)
with the same eigenvalue. In the same way one concludes: \( f \) is an eigenfunction
of the operator
\( \frac{1}{\sqrt{\lambda \exp \beta }}(m_{\cosh (2\sqrt{\beta \pi }x)}\circ \cK^{''}_{c}) \)
with eigenvalue \( \rho  \) iff the function
\( \nu _{\sqrt{\frac{\beta }{\pi }}}\circ Bf \)
is an eigenfunction of the operator
\( \cL _{\beta }:{\cal {H}}L^{2}(\C ,\mu _{\frac{\mid \beta \mid }{\pi }})\rightarrow
{\cal {H}}L^{2}(\C ,\mu _{\frac{\mid \beta \mid }{\pi }}) \)
with eigenvalue \( \rho  \).
\qed

Hence given an eigenfunction \( f=f(x) \) of the operator
\( \frac{1}{\sqrt{\lambda \exp \beta }}(m_{\cosh (2\sqrt{\beta \pi }x)}\circ \cK^{''}_{c}) \)
the corresponding eigenfunction \( F=F(z) \) of the operator
\( \cL _{\beta }:{\cal {H}}L^{2}(\C ,\mu _{\frac{\mid \beta \mid }{\pi }})\rightarrow
{\cal {H}}L^{2}(\C ,\mu _{\frac{\mid \beta \mid }{\pi }}) \)
has the following explicit form:
\[
F(z)=2^{\frac{1}{4}}\int _{\R }f(\xi )\exp \left( 2\sqrt{\pi \beta }\xi z-\pi \xi ^{2}
     -\frac{\beta }{2}z^{2}\right) d\xi .
\]
To relate finally the eigenfunctions \( f \) of the operator \( m_{\cosh (2\sqrt{\beta \pi }x)}\circ \cK^{''}_{c} \)
to those of the Kac-Gutzwiller integral operator \( {\cal {K}}_\beta \) with kernel
\( {\cal {K}}_\beta(\xi ,\eta ) \) we need

\begin{Prop}\label{Eigenwert1}
Let \(H\colon L^{2}(\R,d\xi)\to L^{2}(\R,d\xi)\) be an integral operator
with kernel \(H(\xi,\eta)\) and \(g\colon \R\to ]0,\infty[\) a smooth function.
Then the following statements are equivalent:
\begin{enumerate}
\item[{\rm (1)}] \(f\) is an eigenfunction of the integral operator
\(m_{g^{2}}\circ H\) with kernel \(g^{2}(\xi)H(\xi,\eta)\) for the eigenvalue \(\rho\).
\item[{\rm (2)}] \(\frac{f}{g}\) is an eigenfunction of the integral
operator \(H'\) with kernel \(g(\xi)H(\xi,\eta)g(\eta)\) for the eigenvalue \(\rho\).
\end{enumerate}
\end{Prop}

\pf
This follows from the calculations
\begin{eqnarray*}
\rho f(\xi ) & = & \int _{\R }g^{2}(\xi )H(\xi ,\eta )f(\eta )d\eta \\
 & = & g(\xi )\int _{\R }g(\xi )H'(\xi ,\eta )g(\eta )\frac{f(\eta )}{g(\eta )}d\eta \\
 & = & g(\xi )H({\textstyle\frac{f}{g}})(\xi )
\end{eqnarray*}
 and
\begin{eqnarray*}
\rho \frac{f(\xi )}{g(\xi )} & = & \int _{\R }g(\xi )H(\xi ,\eta )g(\eta )\frac{f(\eta )}{g(\eta )}d\eta \\
 & = & \frac{1}{g(\xi )}\int _{\R }g^{2}(\xi )H(\xi ,\eta )f(\eta )d\eta \\
 & = & \frac{1}{g(\xi )}m_{g^{2}}\circ H(f)(\xi ).
\end{eqnarray*}
 \qed

\begin{Thm}\label{connection}
For real \(\beta\) the Ruelle operator \(\cL_{\beta}\) and the Kac-Gutzwiller operator
\( {\cal {G}_{\beta }} \) with kernel
\((\lambda\exp{\beta})^{-\frac{1}{2}}\cal{K_{\beta}(\xi,\eta)}\) have the same eigenvalues (counted with multiplicities)
and hence the same spectrum. Their eigenfunctions \(F=F(z)\) and \(f=f(\xi)\) to the
eigenvalue \(\rho\) are related as follows:
\begin{enumerate}
\item[{\rm(1)}]
\[
F(z)=(8\pi)^\frac{1}{4}\int_{\R}\sqrt{\cosh{2\sqrt{\pi\beta}}\xi} \ f(2\sqrt{\pi}\xi)
              \exp\left({2\sqrt{\pi\beta}}\xi z-\pi\xi^{2}-\frac{\beta}{2}z^{2}\right) d\xi.
\]
\item[{\rm(2)}]
\[
f(\xi)=\left(\frac{1}{2\pi}\right)^{\frac{1}{4}}\frac{1}{\sqrt{\cosh(2\sqrt{\pi\beta}\xi)}}
\int_{\C} F\left(\sqrt{\frac{\pi}{\beta}}z\right)
\exp\left(\sqrt{\pi} \xi z^{*}-\frac{1}{4}\xi^{2}
-\frac{\pi}{2}{z^{*}}^{2}-\pi z^{2}\right)dz.
\]
 \end{enumerate}
\end{Thm}

\pf
According to Proposition \ref{Ruelle}
 the operators
\( \frac{1}{\sqrt{\lambda \exp \beta }}
(m_{\cosh (2\sqrt{\beta \pi }x)}\circ \cK^{''}_{c}) \) with  $c=2\sqrt \pi$
and \( \cL _{\beta } \) are conjugate via \( \nu _{\sqrt{\frac{\pi }{\beta }}}\circ B \).
Proposition \ref{Eigenwert1} shows that the operator
\( m_{\cosh (2\sqrt{\beta \pi }x)}\circ \cK^{''}_{c} \)
has the same spectrum as the operator
\[
m_{\sqrt{\cosh 2\sqrt{\pi \beta }x}}\circ \cK^{''}_{c}\circ
m_{\sqrt{\cosh 2\sqrt{\pi \beta }x}}
\]
with kernel
\( \sqrt{\cosh 2\sqrt{\pi \beta }x\: \cosh 2\sqrt{\pi \beta }y}\ \cK^{''}_{c}(x,y) \).
But this operator is conjugate to the operator \( {\cal {K}}_\beta \) with kernel
\( {\cal {K}}_\beta(\xi ,\eta ) \) through the map \( R_{c} \) with \( c=2\sqrt{\pi }. \)
Hence, if \( f\in L^{2}(\R ,d\xi ) \) is an eigenfunction of the Kac-Gutzwiller
operator \( \frac{1}{\sqrt{\lambda \exp \beta }}\,{\cal {K}}_{\beta } \) with
eigenvalue \( \rho  \), then
\[
(R_{2\sqrt{\pi }}f)(x)=\sqrt{2\sqrt{\pi }}f(2\sqrt{\pi }x)
\]
is an eigenfunction of the operator
\( \frac{1}{\sqrt{\lambda \exp \beta }}\,
\left(m_{\sqrt{\cosh 2\sqrt{\pi \beta }x}}\circ \cK^{''}_{c}
\circ m_{\sqrt{\cosh 2\sqrt{\pi \beta }x}}\right) \)
with the same eigenvalue. Hence
\( \sqrt{\cosh 2\sqrt{\pi \beta }x}\: (R_{2\sqrt{\pi }}f)(x) \)
is an eigenfunction of the operator
\[ {\textstyle\frac{1}{\sqrt{\lambda \exp \beta }}}\,
\left(m_{\cosh (2\sqrt{\beta \pi }x)}\circ \cK^{''}_{c}\right) \]
with eigenvalue \( \rho  \) by Proposition \ref{Eigenwert1}. But then
\[
F(z)=(8\pi )^{\frac{1}{4}}\int _{\R }\sqrt{\cosh 2\sqrt{\pi \beta }\xi }\
f(2\sqrt{\pi }\xi )
\exp \left(2\sqrt{\pi \beta }\xi z-\pi \xi ^{2}-\frac{\beta }{2}z^{2}\right)\: d\xi
\]
is an eigenfunction of the operator \( \cL _{\beta } \) with eigenvalue \( \rho  \).
On the other hand starting with an eigenfunction of the operator \( \cL _{_{\beta }} \) we
know that \( h(x)=(B^{-1}\circ \nu _{\alpha ^{-1}}F)(x) \) is an eigenfunction
of the operator
\[
{\textstyle\frac{1}{\sqrt{\lambda \exp \beta }}}\
\left(m_{\sqrt{\cosh 2\sqrt{\beta\pi}x}}\circ {\cal K}^{''}_{c}\circ m_{\sqrt{\cosh 2\sqrt{\beta\pi}x}}\right)
\]
which is again conjugate to ${\cal K}_{\beta }$ via the map \( R^{-1}_{c} \)
and hence \( R^{-1}_{c}(\sqrt{\cosh 2\sqrt{\beta\pi}x}^{-1}h)(\xi ) \) 
is an eigenfunction
of the operator \( {\cal G}_{\beta } \). Inserting all the transformations
involved we finally get for the corresponding eigenfunction \( f=f(\xi ) \)
\[
f(\xi )
=({\textstyle\frac{1}{2\pi }})^{\frac{1}{4}}
{\textstyle\frac{1}{\sqrt{\cosh (2\sqrt{\pi \beta })\xi }}}
\int _{\mathbb {C}}F\left({\textstyle\sqrt{\frac{\pi }{\beta }}z}\right)
\exp \left(\sqrt{\pi }\xi z^{*}-\frac{1}{4}\xi ^{2}-\frac{\pi }{2}z^{*2}
            -\pi \mid z\mid ^{2}\right)\: dz.
\]
\qed

\section{Zeros of the Ruelle zeta function}

{}From our discussion it follows that for real \( \beta  \) the eigenvalues of
the Ruelle operator \( \cL _{\beta } \) are real. According to (\ref{Freddet})
the zeros of the
Ruelle zeta function \( \zeta _{R}(\beta ):=\zeta _{R}(1,\beta ) \) are located
at those values of \( \beta  \), where the operator \( \cL  \)\( _{\beta } \)
has eigenvalue \( \rho =\lambda ^{-1} \) and not at the same time the eigenvalue
\( \rho =1 \). In the following we will show for \( 0<\lambda <\frac{1}{2} \)
that there exist infinitely many real \( \beta  \)-values such that \( \cL _{\beta } \)
has the eigenvalue \( \rho =\lambda ^{-1} \) and therefore \( \zeta _{R}(\beta ) \)
presumably has infinitely many ``nontrivial'' zeros on the real line unless unexpected
cancellations take place. For the special case \( \lambda =\frac{1}{2} \) we
can show furthermore that there exists infinitely many values of \( \beta  \)
on the line \( \Re \beta =\ln 2 \) such that the operator \( \cL  \)\( _{\beta } \)
has eigenvalue \( \rho =\lambda ^{-1} \) and hence presumably \( \zeta _{R}(\beta ) \)
has infinitely many trivial zeros on this line \( \Re \beta =\ln 2 \) unless
the aforementioned cancellations take place. We expect this pole and zero structure
for Ruelle's zeta function for the Kac-Baker model to be true indeed for generic
\( 0<\lambda <1 \).

Consider first the case \( \beta =0 \). The spectrum of the operator \( \cL _{0} \)
can be determined explicitly and is given by the numbers \( 2\lambda ^{k},
k=0,1,2,\ldots \).
The corresponding eigenfunctions are polynomials of degree \( k \). Hence the
Ruelle function \( \zeta _{R}(z,\beta ) \) at this point can be calculated
and turns out to be \( (1-2z)^{-1} \). This is just the Artin-Mazur zeta function
for the subshift of finite type over two symbols. The Ruelle function
\( \zeta _{R}(1,\beta ) \)
hence takes the special value \( -1 \) at the point \( \beta =0 \) independently
of the parameter \( \lambda  \). To show now the existence of infinitely
many real \( \beta  \)-values such that the operator \( \cL _{\beta } \) has
eigenvalue \( \rho =\lambda ^{-1} \) it is enough to show that for
\( \beta \rightarrow +\infty  \),
respectively for \( \beta \rightarrow -\infty  \), infinitely many eigenvalues
become larger than \( \lambda ^{-1} \), since the eigenvalues depend analytically
on \( \beta  \).
The asymptotic behaviour of the eigenvalues has been
determined by B.\ Moritz in his unpublished diploma thesis \cite{Mo89}.
To formulate his result we have to introduce the parity operator
\( P:B(D)\rightarrow B(D) \)
defined via
\[
Pf(z)=f(-z).\]
It is straightforward to see that all eigenspaces of the operator
\( \cL  \)\( _{\beta } \)
have a basis consisting of eigenfunctions being either even (eigenvalue \( +1 \))
or odd (eigenvalue \( -1 \)) under the parity operator \( P \). But then the
spectrum of this operator is just the union of the spectra of the two following
operators \( \cL ^{+}_{\beta }:B(D)\rightarrow B(D) \) and
\( \cL ^{-}_{\beta }:B(D)\rightarrow B(D) \)
defined by
\[
\cL ^{+}_{\beta }f(z)
:=\exp \left(\beta z\right)\: f(\lambda +\lambda z)
 +\exp\left( -\beta z\right)\: f(\lambda -\lambda z),
\]
respectively
\[
\cL ^{-}_{\beta }f(z)
:=\exp \left(\beta z\right)\: f(\lambda +\lambda z)
 -\exp \left(-\beta z\right)\: f(\lambda -\lambda z).
\]
Obviously the eigenfunctions of these two operators are even, respectively odd.
Their eigenvalues will be called even, respectively odd. The main result of Moritz
is then the following

\begin{Prop}\label{Moritz}
For arbitrary N consider the N even and odd eigenvalues \(\rho _{i}\) of the operator
\({\cal L}_{\beta}\) largest
in absolute value. For \(0<\lambda<\frac{1}{2}\) these eigenvalues behave for
\(\beta\to\infty\) as \(\lambda^{i} \exp\left({\frac{\lambda }{1-\lambda}\beta}\right)\).
For \(\beta\to-\infty \) the N even eigenvalues largest
in absolute value behave like
\((-1)^{i}\lambda^{i}\exp\left({-\frac{\lambda }{1+\lambda}\beta}\right)\),
whereas the odd eigenvalues behave like
\((-1)^{i+1}\lambda^{i}\exp\left(-{\frac{\lambda }{1+\lambda}\beta}\right) \).
\end{Prop}

\pf
Since the arguments for the odd eigenvalues \( \rho ^{-} \) are similar to
the even eigenvalues \( \rho ^{+} \) we restrict ourselves to the even eigenvalues.
Consider  the case \( \beta \rightarrow \infty  \) first. If we write in this
case the eigenfunction \( f \) with eigenvalue \( \rho  \) for the operator
\( L^{+}_{\beta } \) in the form
\( f(z)=\exp\left( \frac{\beta z}{1-\lambda }\right)\: u(z) \),
one finds with
\( \overline{\rho }=\exp\left( -\frac{\beta \lambda }{1-\lambda }\right)\: \rho  \)
and \( \overline{\beta }=\frac{2\beta }{1-\lambda } \) by a simple calculation
\[
\overline{\rho }u(\lambda +z')
=u(\lambda +\lambda ^{2}+\lambda z')
 +\exp\left( -\overline{\beta }\lambda\right)\exp\left( -\overline{\beta }z'\right)\:
 u(\lambda -\lambda ^{2}-\lambda z'),
\]
where we have also replaced the argument \( z \) by \( z'+\lambda  \). For
the function \( h(z):=u(\lambda +z) \) one then finds the equation
\[
h(\lambda z+\lambda ^{2})
+\exp\left( -\overline{\beta }\lambda\right) \exp\left(-\overline{\beta }z\right)\:
h(-\lambda ^{2}-\lambda z)
=\overline{\rho }h(z).
\]
Hence the function \( h \) is an eigenfunction with eigenvalue \( \overline{\rho } \)
of the operator \( {\cal {T}}_{\overline{\beta }}:B(D_{R})\rightarrow B(D_{R}) \)
defined as
\[
{\cal {T}}_{\overline{\beta }}h(z)
=h(\lambda z+\lambda ^{2})
 +\exp\left(-\overline{\beta }\lambda\right) \exp\left(-\overline{\beta }z\right)\:
 h(-\lambda ^{2}-\lambda z),
\]
where \( D_{R}=\{z\in \mathbb {C}\mid \ |z| <R\} \) and
\( R>\frac{\lambda ^{2}}{1-\lambda } \).
This operator is nuclear and its spectrum is closely related to the
one of the operator \( {\cal L}^{+}_{\beta } \). Writing now
\[
{\cal {T}}_{\overline{\beta }}
={\cal {T}}^{+}_{\overline{\beta }}+{\cal {T}}^{-}_{\overline{\beta }}
\]
with the obvious definitions of the two operators \( {\cal {T}}^{+}_{\overline{\beta }} \)
and \( {\cal {T}}^{-}_{\overline{\beta }} \) one finds for the norm of the
operator \( \cal {T}^{-}_{\overline{\beta }} \)
\[
\parallel {\cal {T}}^{-}_{\overline{\beta }}\parallel
\leq \exp\left( -\overline{\beta }\lambda\right) \exp\left( \overline{\beta }R\right).
\]
Since for \( \lambda <\frac{1}{2} \) one can choose
\( \lambda >R>\frac{\lambda ^{2}}{1-\lambda } \)
we find
\( \lim _{\overline{\beta }\rightarrow \infty }\parallel
{\cal {T}}^{-}_{\overline{\beta }}\parallel =0. \)
Hence the operator \( {\cal {T}}_{\overline{\beta }} \) approaches for
\( \overline{\beta }\rightarrow \infty  \)
in norm the operator \( {\cal {T}}^{+}_{\overline{\beta }} \) and hence also
their spectra are identical in this limit. The spectrum of the operator
\( {\cal {T}}^{+}_{\overline{\beta }} \),
however, is given by the numbers \( \{\lambda ^{i}\mid i=0,1,\ldots\} \).
This shows that for large \( \beta  \) the even eigenvalues \( \rho ^{+} \)
of the operator \( {\cal L}_{\beta } \) behave like
\( \lambda ^{i}\exp \left(\frac{\lambda \beta }{1-\lambda }\right) \).
To find the behaviour of the even eigenvalues of the operator \( \cL _{\beta } \)
for large negative values of \( \beta  \) consider the operator
\( \widetilde{\cL }^{+}_{\beta } \)
defined as
\[
\widetilde{\cL }^{+}_{\beta }f(z)
:=\exp\left( -\beta z\right)\: f(\lambda +\lambda z)
  +\exp\left( \beta z\right)\: f(\lambda -\lambda z)
\]
and the behaviour of its eigenvalues for large positive values of \( \beta  \).
Writing its eigenfunction \( f \) with eigenvalue \( \rho  \) as
\( f(z)=\exp\left( \frac{\beta z}{1+\lambda }\right)u(z) \)
one finds in this case with
\( \overline{\rho }:=\exp\left( -\frac{\beta \lambda }{1+\lambda }\right) \)
and
\( \overline{\beta }:=\frac{2\beta }{1+\lambda } \) the equation
\[
h(-\lambda ^{2}-\lambda z)
+\exp\left( -\overline{\beta }\lambda\right) \exp\left( -\overline{\beta }z\right)\:
h(\lambda ^{2}+\lambda z)
=\overline{\rho }h(z),
\]
where we have introduced again the function \( h(z)=u(\lambda +z) \). An argument
completely analogous to the former case then shows that the large \( \overline{\beta } \)
behaviour of the eigenvalues \( \overline{\rho } \) is determined by the operator
\( {\cal {\widetilde{T}}}^{+}_{\overline{\beta }}h(z)=h(-\lambda ^{2}-\lambda z) \).
Its spectrum is given by the numbers \( \{(-\lambda )^{i}\mid i=0,1,\ldots\} \).
Hence the even eigenvalues of the operator \( \cL _{\beta } \) behave for large
negative \( \beta  \) as
\( (-\lambda )^{i}\exp\left( -\frac{\beta \lambda }{1+\lambda } \right)\).
This proves the claim.
\qed

Let us calculate the trace of the Ruelle operator \( \cL _{\beta } \) for large
positive or negative values of \( \beta  \). Adding up \( N \) of the asymptotic
eigenvalues one finds for large positive \( \beta  \):
\[ \trace\cL _{\beta }
 \underset{\beta\to\infty}{\sim}
  \frac{1-\lambda ^{N+1}}{1-\lambda }\ \exp\left( \frac{\beta \lambda }{1-\lambda }\right).
\]
For large negative \( \beta  \) on the other hand one finds when adding up
\( N \) asymptotic eigenvalues
\[ \trace\cL _{\beta }
\underset{\beta\to-\infty}{\sim}
O(\lambda ^{N+1}).
\]
{}From this we expect that the result of Moritz on the asymptotic behaviour of
the eigenvalues indeed is true for general \( 0<\lambda <1 \).

Now we can prove the following proposition which determines the location of
zeros and poles of the Ruelle zeta function for the Kac model:

 \begin{Prop}\label{zeros}
For any \(\lambda\) with \(0<\lambda<1\) the Fredholm determinant
\(\det(1-\lambda\cal{L}_\beta)\) has infinitely many zeros on the real line \(\R\).
For \(\lambda=\frac{1}{2}\) there are infinitely many zeros on the line
\(\Re\lambda=\ln 2\).
\end{Prop}

\pf
Since for real \( \beta  \) all the eigenvalues \( \rho (\beta ) \) of the
operator \( \cL _{\beta } \) are real and analytic in \( \beta  \), it follows
from the asymptotic behaviour of the eigenvalues that infinitely many of them
must take the value \( \lambda ^{-1} \) for positive and negative values of
the parameter \( \beta  \). On the other hand an easy calculation shows that
for the special value \( \lambda =\frac{1}{2} \) the function \( f(z)=\sinh(2\beta z) \)
is an eigenfunction of the operator \( \cL _{\beta } \) for all \( \beta \in \mathbb {C} \)
with eigenvalue \( \rho =\exp \beta  \). But this eigenvalue takes the value
\( \lambda ^{-1}=2 \) just for \( \beta =\beta _{n}=\ln 2\: +2\pi in \).
\qed

The zeros of the Fredholm determinant \( \det (1-\cL _{\beta }) \) on the real
line are certainly ``nontrivial'', whereas the ones on the line \( \Re\beta =\ln 2 \)
which even are equidistant could be called ``trivial'' ones. Since the Ruelle
zeta function \( \zeta _{R}(z) \) has the representation
\( \zeta _{R}(z)=\frac{\det (1-\lambda \cL _{\beta })}{\det(1-\cL _{\beta })} \)
this function has infinitely many nontrivial zeros and poles on the real line,
whereas for the special value \( \lambda =\frac{1}{2} \) infinitely many trivial
zeros lie on the line \( \Re z=\ln 2 \). Obviously accidental cancellations with
the zeros of \( \det (1-\cL _{\beta }) \) which determine the poles of the
Ruelle zeta function could destroy some of these zeros.

Since the location of the nontrivial zeros of the Ruelle zeta function is not
known explicitly one has to determine them numerically. Indeed in his paper
\cite{Gu82} M.\ Gutzwiller derived explicit formulas for the matrix
elements of more or less the operator 
\( \frac{1}{\sqrt{\lambda \exp \beta }}\ 
m_{\cosh (2\sqrt{\beta \pi }x)}\circ {\cal {K}}^{''}_{c} \)
in the basis given by the Hermite functions \( h_{k}(\xi ) \) in the space
\( L^{2}({\mathbb {R}},d\xi ) \). After a conjugation this matrix becomes a
symmetric matrix \( {\cal {B}}_{\beta } \) whose matrix elements 
\( {\cal {B}}_{n,m}(\beta ) \) have
the following form:
\[
{\cal {B}}_{n,m}(\beta )=2(n!m!)^{-\frac{1}{2}}
\exp\left( -\frac{(n+m)\gamma }{2}\right){\textstyle\frac{M!}{(2\mu )!}}
\Phi (2\mu -M,2\mu +1;-\beta ),
\]
where \( \Phi  \) denotes the confluent hypergeometric function,
\( M=\max\{m,n\} \), and the number \( 2\mu =|m-n|  \) must be
even. For \(|m-n|\) odd the matrix elements \( {\cal {B}}_{n,m}(\beta ) \)
vanish.  Since the number \( 2\mu -M \) is a non-positive integer the above confluent 
hypergeometric
function is just proportional to the Laguerre polynomial according to the formula
\[
\Phi (-n,a+1;x)=\left( _{n}^{n+a}\right) ^{-1}L^{a}_{n}(x).
\]
Inserting this relation into the expression for the matrix elements \( {\cal {B}}_{n,m} \)
one finally gets
\[
{\cal {B}}_{n,m}(\beta )
=2(n!m!)^{-\frac{1}{2}}\exp\left( -\frac{(n+m)\gamma }{2}\right)
 {\textstyle \frac{M!}{(2\mu !)}}
 \left( _{M-2\mu }^{M}\right) ^{-1}L_{M-2\mu }^{2\mu }\left(-\beta\right).
\]
This allows one to calculate the traces of the iterates of the Kac-Gutzwiller
Operator as the sum over the diagonal elements of this matrix. For instance
for the trace of \( {\cal {G}}_{\beta } \) one obtains in this way
\[
\trace\: {\cal {G}}_{\beta }
=2\sum ^{\infty }_{m=0}\exp\left( -\gamma m\right)\: L_{m}^{0}(-\beta )
=2\frac{1}{1-\lambda }\exp\left( \frac{\beta \lambda }{1-\lambda }\right),\]
which is just the definition of the generating function for the Laguerre polynomials
\( L_{m}(-\beta )=L^{0}_{m}(-\beta ). \) In complete analogy one can relate
the traces of the iterates of the Kac-Gutzwiller operator \( {\cal {G}}_{\beta } \)
calculated via the matrix \( {\cal {B}}(\beta ) \) to the expressions calculated
via the Ruelle operator \( {\cal {L}}_{\beta } \) and its iterates. One arrives
thereby at rather complicated formulas which can be interpreted as generating
functions for powers of the Laguerre polynomials \( L_{M-2\mu }^{2\mu }(-\beta ) \). 
Whether these formulas are known in the literature for the Laguerre polynomials
is not known to us. 

A detailed numerical study of the zeros and poles of the Ruelle zeta function
by the above matrices is under way.

{\bf Acknowledgements}: This work has been supported by the Deutsche Forschungsgemeinschaft through the DFG Forschergruppe ``Zeta Functions and Locally Symmetric Spaces''. D. M. thanks the MPI for Mathematics in Bonn for the kind hospitality extended to him. There the final version of this paper has been prepared.


\begin{thebibliography}{Mo89}

\bibitem[AtBo67]{AtBo67} M. Atiyah, R. Bott. A Lefschetz fixed point formula for elliptic
                         complexes I. Ann. Math. \( \underline{86} \), 374--407 (1967).

\bibitem[Ba61]{Ba61} G. Baker. One dimensional order-disorder model which approaches a second order phase transition. 
                     Phys. Rev. \( \underline{122} \),
                     1477--1484 (1961).

\bibitem[Be86]{Be86} M. Berry. Riemann's zeta function, a model of quantum chaos. Lect.
                     Notes in Physics \( \underline{263} \), Springer Verlag, Berlin, 1986.

\bibitem[Co96]{Co96} A. Connes. Formule de trace en geometrie non commutative et hypothese
                     de Riemann. C.R. Acad. Sci. Paris Ser I
                     \( \underline{323} \), 1231--1236 (1996).

\bibitem[Cr46]{Cr46} H. Cram\'er. \textit{Mathematical Methods in Statistics.}  Princeton
                     Univ. Press, 1946.

\bibitem[De99]{De99} C. Deninger. On dynamical systems and their possible significance
                    for arithmetic geometry.
                    Progr. Math. \( \underline{171} \), 29--87 (1999).

\bibitem[Fo84]{Fo84}G.B. Folland. \textit{Real Analysis.} Wiley, New York, 1984.

\bibitem[Fo89]{Fo89}G.B. Folland. \textit{Harmonic Analysis in Phase Space.}
                    Princeton University Press, 1989.

\bibitem[Gu82]{Gu82} M. Gutzwiller. The quantization of a classically ergodic system.
                      Physica \( \underline{5D} \), 183--207 (1982).

\bibitem[Is01]{Is01} S.Isola. On the spectrum of Farey and Gauss maps. Preprint Univerity
                      of Camerino, 2001.

\bibitem[Ka59]{Ka59} M. Kac. On the partition function of a one-dimensional lattice gas.             Phys. Fluids \( \underline{2} \),
                     8--12 (1959).

\bibitem[Ka66]{Ka66} M. Kac. The mathematical mechanism of phase transitions.
                     In \textit{Brandeis University Summer Institute in Theoretical Physics},
                     vol 1, 245--305 eds.: H. Chretien et al. Gordon \& Breach, New York,
                     1966.

\bibitem[Ma80]{Ma80} D. Mayer. The Ruelle-Araki transfer operator in classical statistical
                     mechanics. Lect. Notes in Physics \( \underline{123} \)(1),
                     Springer Verlag, Berlin, 1980.

\bibitem[Ma90]{Ma90} D.Mayer. On the thermodynamic formalism for the Gauss map. Commun.
                       Math. Phys. \( \underline{130} \), 311--333 (1990).

\bibitem[Ma91]{Ma91} D. Mayer. The thermodynamic formalism approach to Selberg's zeta
                       function for \( PSL(2,\mathbb {Z}) \).
                       Bull.\ Am.\ Math.\ Soc. \( \underline{25} \), 55--60 (1991).

\bibitem[Mo89]{Mo89}B. Moritz. \textit{Die Transferoperator--Methode in der
                    Behandlung des Kacschen Spinmodells.} Diplomarbeit, RWTH Aachen, 1989.

\bibitem[PTMT94]{PTMT94} A. Pentek, Z. Toroczkai, D. Mayer, T. Tel. The Kac model from
                       a dynamical system's point of view.
                       Phys. Rev. E \( \underline{49} \), 2026--2040 (1994).

\bibitem[Ro86]{Ro86} P. Robba. Une introduction naive aux cohomologies de Dwork. Mem.
                       Soc. Math. France \( \underline{23} \), 61--105 (1986).

\bibitem[Ru68]{Ru68} D. Ruelle. Statistical mechanics of a one dimensional lattice gas.
                       Commun. Math. Phys. \( \underline{9} \), 267--278 (1968).

\bibitem[Ru92]{Ru92} D. Ruelle. Dynamical zeta functions: where do they come from and
                       what are they good for? In \textit{Proc. Internat. Congress of
                       Math. Phys. X (Leipzig 1991)}, 43--51, Springer Verlag, Berlin, 1992.

\bibitem[Rug94]{Rug94} H.H. Rugh. On the asymptotic form and the reality of spectra of
                        Perron-Frobenius operators.
                        Nonlinearity \( \underline{7} \), 1055--1066 (1994).
\bibitem[ViMa77]{ViMa77} K. Vishwanathan, D. Mayer. Statistical mechanics of one-dimensional Ising and Potts models with exponential interactions.
                        Physica \(\underline {89A}\), 97--112 (1977).  








\end{thebibliography}
\end{document}